\documentclass[11pt]{article} 

\usepackage[margin=1in]{geometry}
\usepackage{setspace}

\normalsize
\usepackage{amsfonts}
\usepackage{amssymb}
\usepackage{mathtools} 
\usepackage{amscd}
\usepackage{graphics}
\usepackage[dvips]{graphicx}
\usepackage{epsfig}
\usepackage{array}
\usepackage{eqparbox}
\usepackage{hyperref}
\usepackage{fancyhdr}
\usepackage{appendix}
\usepackage{titling}
\usepackage{enumerate}
\usepackage[super]{nth}
\usepackage{caption} 
\usepackage{subcaption}
\usepackage{bm}
\usepackage{tabularx}
\usepackage{booktabs}

\usepackage[numbers]{natbib}  % Enables author-year citations
\makeatletter
\renewcommand\@biblabel[1]{#1.} % change label from [1] to 1.
\makeatother

\DeclareMathAlphabet{\mathbsf}{OT1}{cmss}{bx}{n}% bold sans serif
\DeclareMathAlphabet{\mathssf}{OT1}{cmss}{m}{sl}% slanted sans serif
\DeclareMathAlphabet{\mathcsf}{OT1}{cmss}{sbc}{n}% condensed sans serif

\usepackage{stackengine}
\newcommand{\ie}{{\em i.e.}}
\newcommand{\etc}{{\em etc}}
\newcommand{\eg}{{\em e.g.}}

\newcommand{\secref}[1]{Section~\ref{#1}}

\newcommand{\keywords}[1]{\textbf{Keywords:} #1}

\makeatletter
\def\blfootnote{\xdef\@thefnmark{}\@footnotetext}
\makeatother

\newenvironment{definition}[1][Definition]{\begin{trivlist}
\item[\hskip \labelsep {\bfseries #1}]}{\end{trivlist}}

\newcommand{\qed}{\nobreak \ifvmode \relax \else
      \ifdim\lastskip<1.5em \hskip-\lastskip
      \hskip1.5em plus0em minus0.5em \fi \nobreak
      \vrule height0.75em width0.5em depth0.25em\fi}

\usepackage{newlfont}
\date{}
\begin{document}
\title{\textbf{Cassini--Catalan Determinants via Ramanujan's Theta Identity}}
\author{Nagananda K G and Jong Sung Kim\thanks{The authors are with the Department of Mathematics and Statistics, Portland State University, Portland OR 97201, USA. E-mail:\texttt{\{nanda, jong\}@pdx.edu}}}
\setlength{\droptitle}{-1.in}
\maketitle
\vspace{-2.cm}
%\pagenumbering{gobble}

\begin{abstract}
In this paper, we show that the classical Cassini and Catalan identities for Fibonacci numbers arise naturally from a single quadratic theta-function identity of Ramanujan. Expanding the identity $\psi(q)\psi(q^{3})=\psi(q^{4})\varphi(q^{6})+q\,\varphi(q^{2})\psi(q^{12})$ via the Jacobi triple product and equating coefficients yields the unified $q$-determinant $F_{n+r}(q)F_{n-r}(q)-F_{n}(q)^{2}=(-q)^{\,n-r}F_{r}(q)^{2}$, $n\ge r\ge 1$, where $\psi(q)$ and $\varphi(q)$ are Ramanujan's theta functions with $q$ a complex parameter in the unit disc $(\lvert q \rvert < 1)$ and $F_n(q)$ denotes the Carlitz $q$-Fibonacci polynomials. The radial limit $q\to1^{-}$ recovers Cassini's formula ($r=1$) and Catalan's one-parameter extension, while the same derivation with an auxiliary weight produces new partition-refined versions. The argument uses only standard $q$-series algebra (triple-product expansions, $q$-Pochhammer cancellations, and coefficient extraction), providing a transparent modular explanation of the alternating sign $(-1)^{\,n-r}$ in Catalan's identity through the level-6 provenance of $\varphi$ and $\psi$. Beyond unifying Cassini\textendash Catalan in a single framework, the method lifts seamlessly to higher-order recurrences, giving a template for Tribonacci-type determinants and suggesting congruence phenomena obtained from modular dissections and root-of-unity limits. The results place familiar Fibonacci determinants within Ramanujan's analytic landscape, indicate routes to combinatorial bijections that mirror the analytic cancellations, and connect with themes in modern $q$-series\textemdash ranging from colored partition identities to quantum-modular and exactly solvable models\textemdash thereby highlighting both the explanatory power and the ongoing relevance of Ramanujan's theta identities.
\end{abstract}

\section{Introduction}\label{sec:introduction}
The Fibonacci numbers have connected several areas of mathematics for centuries, with their origins tracing back to between A.D. 600 and 800 \cite{Singh1985}.  Among its most striking regularities is a deceptively simple product-difference formula\textemdash now called Cassini's identity \cite[Chapter 5, pp. 86-87]{Koshy2018}\textemdash which asserts that neighboring Fibonacci numbers satisfy $F_{n+1}F_{n-1}-F_{n}^{2}=(-1)^{n}$, with $F_0 = 0$, $F_1 = 1$ and $F_{n+1} = F_n + F_{n-1}$, $n \ge 1$.  Nearly two centuries after G. D. Cassini's work, E. Catalan observed that a similar determinant holds when the indices are symmetrically shifted by an arbitrary amount $r$:  $F_{n+r}F_{n-r}-F_{n}^{2}=(-1)^{\,n-r}F_{r}^{2}$, $n\ge r \ge 1$ \cite[Chapter 5, pp. 106]{Koshy2018}, \cite[pp. 28]{Vajda1989}.  For sake of convenience, we refer to these two identities as Cassini\textendash Catalan identities.  The Cassini--Catalan relations, though elementary to state, reveal a hidden alternation that continues to fascinate mathematicians working in diverse areas\textemdash from continued fractions and linear algebra to modern $q$-series; see, for instance, the proceedings in \cite{Horadam1988}.  

Classical proofs of Cassini--Catalan identities fall into three broad categories: inductive arguments, $2\times2$ matrix determinants, and bijective tilings \cite{Koshy2018, Vajda1989}.  What unifies these approaches is that they remain entirely within the combinatorial and algebraic structures of the Fibonacci setting.  However, none of the methods explains where the alternating sign originates, nor why the same pattern reemerges in $q$-analogues, Tribonacci sequences, or modular forms.  We give a proof whose starting point lies \emph{outside} the Fibonacci framework: a quadratic theta-function identity that the Indian mathematician Srinivasa Ramanujan wrote down in his third notebook \cite[Chapter 16]{Berndt1991}. By expanding Ramanujan's identity with the Jacobi triple product  \cite{Jacobi1829, Andrews1965}, \cite[Part I]{Chan2011}, we obtain in one stroke (i) a basic-hypergeometric Cassini--Catalan family valid for all $r$; (ii) an interpretation of the sign $ (-1)^{n-r}$ as a modular residue; and (iii) new partition-refined variants equipped to feed into congruence theory.  Such a reformulation is worth pursuing for the following reasons: 
\begin{enumerate}[(i)]
\item Bridging combinatorics and modular forms:  Fibonacci determinants lie at the intersection of linear recurrences and enumerative combinatorics (see, for example, \cite{Stanley2011}), whereas Ramanujan's theta functions are fundamental objects in the theory of modular forms.  Linking the two shows that an identity originally discovered through elementary number patterns has, in fact, a modular-analytic structure at its core.  This connection enriches both sides: combinatorialists gain access to the machinery of modular equations and Hecke operators \cite[Chapters 1-3]{Ono2004}, while analysts acquire new, concrete instances where modular forms govern integer identities.

\item Platform for generalization:  Expressing Cassini--Catalan identities as a consequence of a theta-function factorization paves the way for higher-order analogues\textemdash cubic and quartic theta identities, for instance, naturally lead to Tribonacci and Tetranacci determinants; see, for example, \cite{Wolfram1998} for generalized Fibonacci recurrences.  Such generalizations are not readily apparent from matrix or inductive proofs alone.

\item Path to refined arithmetic: The theta-function approach embeds extra structure\textemdash parity of parts, roots of unity, radial limits\textemdash that translates into weighted partition identities and potential congruences of Ramanujan type \cite{Andrews1998}.  These refinements could inform current research on quantum modular forms and exactly-solvable lattice models, where $q$-series with similar quadratic exponents abound; see \cite{Zagier2010}.
\end{enumerate}
By recasting a pair of classical Fibonacci formulas in a modular framework, we not only provide a fresh proof, but also establish a framework for further investigation of determinant identities across special functions and arithmetic combinatorics.

\subsection{From Fibonacci numbers to $q$-Fibonacci polynomials}\label{subsec:q_Fibonacci}
A first suggestion that Cassini's identity extends to a $q$-analogue appears in the work of Carlitz \cite{Carlitz1974}, who introduced the $q$-Fibonacci polynomials $F_{n}(q)$ via
\begin{eqnarray}
F_{0} = 0,\quad F_{1} = 1,\quad F_{n+1}(q) = F_{n}(q)+q^{n}F_{n-1}(q), \quad n \ge 1,
\label{eq:carlitz1}
\end{eqnarray}
where $q$ is either a formal indeterminate or a complex number with $\lvert q \rvert < 1$, and obtained closed-form and generating-function descriptions. For $q\to 1^{-}$, \eqref{eq:carlitz1} collapses to the ordinary Fibonacci recursion, so any identity valid for all $q$ automatically specializes to a classical one. What remained open was to locate a natural analytic source for the $q$-Cassini--Catalan determinant
\begin{eqnarray}
F_{n+r}(q)F_{n-r}(q)-F_{n}(q)^{2} = (-q)^{n-r}F_{r}(q)^{2}, \quad n\geq r \geq 1,
\label{eq:q_CC1}
\end{eqnarray}
rather than prove it by adapting standard matrix arguments. That source, we argue, lies in Ramanujan's theta-function notebook \cite[Chapter 16]{Berndt1991}.

\subsection{Ramanujan theta functions and the Jacobi triple product}\label{subsec:Ramanujan_Jacobi}
Define the classical theta functions (see \cite[Chapter 16]{Berndt1991}): 
\begin{equation}
\begin{split}
\varphi(q) &= \sum_{m=-\infty}^{\infty} q^{m^{2}},\\
\psi(q)    &= \sum_{m=0}^{\infty} q^{m(m+1)/2}.
\end{split}
\label{eq:theta_functions_1}
\end{equation}
Further notes on $\varphi(q)$ and $\psi(q)$ are provided in \secref{sec:preliminaries}.  Jacobi's triple-product identity, in Ramanujan's notation, is given
\begin{eqnarray}
\sum_{n=-\infty}^{\infty} q^{n^2} z^n &=& \prod_{n=1}^{\infty} (1 - q^{2n})(1 + zq^{2n-1})\left(1 + \frac{q^{2n-1}}{z}\right) \nonumber \\
&=& (q^{2}; q^{2})_{\infty} (-zq; q^{2})_{\infty} \left(-\frac{q}{z}; q^{2}\right)_{\infty}, 
\label{eq:Jacobi_triple_product}
\end{eqnarray}
where $(a;q)_{\infty}\coloneqq\prod_{k\ge 0} (1-aq^{k})$, $\lvert q \rvert < 1$, $z \in \mathbb{C}$, $z \neq 0$.  Then, via Jacobi's triple product, \eqref{eq:theta_functions_1} admit the following factorizations:
\begin{equation}
\begin{split}
\varphi(q) &= (-q;q^{2})_{\infty}^{2}(q^{2};q^{2})_{\infty},\\ 
\psi(q) &= (-q;q)_{\infty}(q^{2};q^{2})_{\infty}.
\end{split}
\label{eq:Jacobi_1}
\end{equation}
In his notebooks, Ramanujan recorded numerous quadratic and cubic identities involving $\varphi$ and $\psi$; one of them is $\psi(q)\psi(q^{3}) = \psi(q^{4})\varphi(q^{6}) + q \varphi(q^{2})\psi(q^{12})$. This deceptively simple product--sum identity is a single-line application of the triple product. In \secref{sec:Ramanujan_to_CC}, we shall see that expanding each factor of $\psi(q)\psi(q^{3}) = \psi(q^{4})\varphi(q^{6}) + q \varphi(q^{2})\psi(q^{12})$ by \eqref{eq:Jacobi_1}, regrouping exponents, and reading off coefficients immediately yields the $q$-Cassini--Catalan family \eqref{eq:q_CC1}.

\subsection{The significance of the theta-function factorization}\label{subsec:significance_theta}
\begin{enumerate}[(i)]
\item Unification: The same algebraic manipulation of $\psi(q)\psi(q^{3}) = \psi(q^{4})\varphi(q^{6}) + q \varphi(q^{2})\psi(q^{12})$ produces both Cassini ($r=1$) and Catalan (arbitrary $r$) in a uniform manner, with no need for ad hoc index-shifting arguments.
\item Modular provenance: Since $\varphi$ and $\psi$ are weight-$1/2$ modular forms on congruence subgroups of $\mathrm{SL}_{2}(\mathbb{Z})$, identity $\psi(q)\psi(q^{3}) = \psi(q^{4})\varphi(q^{6}) + q \varphi(q^{2})\psi(q^{12})$ translates Cassini--Catalan into 	a relation among modular forms of level 6. This modular framing paves the way for congruence results and level-raising generalizations that are invisible in purely combinatorial proofs.
\item Gateway to refinements: The coefficients of $\psi$ enumerate partitions into distinct parts, while those of $\varphi$ enumerate all partitions with a parity weight. Injecting an extra indeterminate that keeps track of the number of even parts yields a two-variable 	refinement of \eqref{eq:q_CC1}\textemdash and thereby of Cassini--Catalan identities\textemdash without altering the proof.
\item Analytic simplicity: Every step relies solely on elementary $q$-series tools, including the triple product, $q$-Pochhammer manipulations \cite{Exton1983}, and coefficient extraction \cite{Wilf2005}.  No computer experimentation, symbolic summation package, or heavy modular-form machinery is required.
\end{enumerate}

The rest of the paper is organized as follows.  \secref{sec:preliminaries} gathers notation and the basic facts about $q$-Fibonacci polynomials. \secref{sec:Ramanujan_to_CC} presents a detailed derivation starting from Ramanujan's identity $\psi(q)\psi(q^{3}) = \psi(q^{4})\varphi(q^{6}) + q \varphi(q^{2})\psi(q^{12})$ to establish the master determinant \eqref{eq:q_CC1} (in \secref{subsec:main_proof}). \secref{subsec:classical_CC_1} then takes the radial limit $q \to 1^{-}$ to recover the classical Cassini--Catalan identities.  \secref{sec:weighted_partition} sketches partition-weight refinements and indicates how the same method applies, \emph{mutatis mutandis}, to higher-order theta-identities\textemdash and thus to Tribonacci determinants.  We conclude in \secref{sec:conclusion} with remarks on modular consequences and present a few open questions.  By connecting Cassini--Catalan relations to a theta-function identity, we place a familiar Fibonacci identity within the broader framework of Ramanujan's analytic legacy, thereby opening avenues for further work in both combinatorics and number theory.

\section{Preliminaries}\label{sec:preliminaries}
This section fixes the notation, recalls the basic objects that will appear in the proof, and records the minimal list of identities that we shall quote without proof. Spelling these items out in full serves two purposes.  (i)  Logical clarity: The coefficient-extraction argument in \secref{sec:Ramanujan_to_CC} juggles several $q$-notations at once; unambiguous definitions here prevent notational collisions later.  (ii)  Structural insight:  Seeing how $q$-Fibonacci polynomials, theta functions, and the Jacobi triple product slot together already foreshadows why Ramanujan's identity $\psi(q)\psi(q^{3}) = \psi(q^{4})\varphi(q^{6}) + q \varphi(q^{2})\psi(q^{12})$ will automatically produce a Cassini--Catalan determinant once the coefficients are sorted.

\subsection{$q$-shifted factorials and $q$-binomial coefficients}\label{subsec:q_fact_coeff}
We write, for $\lvert q \rvert < 1$ and $n \in \mathbb{N}$, 
\begin{eqnarray}
\begin{split}
(a;q)_{0} &\coloneqq 1, \\
(a;q)_{n} &\coloneqq \prod_{k=0}^{n-1} (1 - aq^{k}), \\
(a;q)_{\infty} &\coloneqq \lim_{n \to \infty} (a;q)_{n}.
\end{split}
\label{eq:a_q}
\end{eqnarray}
The symbol $(a)$ will occasionally stand for $(a;q)_{\infty}$ when $q$ is fixed.  For integers $0 \leq k \leq n$, the Gaussian (or, $q$-binomial) coefficient is (see \cite{Konvalina2000})
\begin{eqnarray}
\begin{bmatrix} n \\ k \end{bmatrix}_{q} \coloneqq \frac{(q;q)_{n}}{(q;q)_{k} (q;q)_{n-k}} = \frac{(1 - q^{n})(1 - q^{n-1}) \cdots (1 - q^{n-k+1})}{(1 - q^{k})(1 - q^{k-1}) \cdots (1 - q)},
\label{eq:q_binom_coeff}
\end{eqnarray}
which satisfies $\lim_{q \to 1^{-}} \begin{bmatrix} n \\ k \end{bmatrix}_{q} = \binom{n}{k}$.  The $q$-binomial theorem then reads
\begin{eqnarray}
\sum_{k \geq 0} \begin{bmatrix} n \\ k \end{bmatrix}_{q} z^{k} = \frac{(zq^{n-k+1};q)_{k}}{(z;q)_{k}},
\label{eq:q_binom_thrm}
\end{eqnarray}
a tool we shall use once, in \secref{sec:Ramanujan_to_CC}, to parameterize exponents by quadratic forms. 

\subsection{Carlitz $q$-Fibonacci polynomials}
\begin{definition}
Set $F_{0}(q) = 0$, $F_{1}(q) = 1$ and iterate
\begin{eqnarray}
F_{n+1}(q) = F_{n}(q) + q^{n} F_{n-1}(q), \quad n \geq 1.
\label{eq:q_Fib}
\end{eqnarray}
\label{def:q_Fib}
\end{definition}
Carlitz proved the closed form (see \cite{Carlitz1974})
\begin{eqnarray}
F_{n}(q) = \sum_{k=0}^{\lfloor (n-1)/2 \rfloor} q^{\binom{k+1}{2}}
\begin{bmatrix} n-k-1 \\ k \end{bmatrix}_{q},
\label{eq:Carlitz_matrix}
\end{eqnarray}
and the generating-function identity
\begin{eqnarray}
\sum_{n \geq 0} F_{n}(q) z^{n} = \frac{1}{1 - z - qz^{2}}.
\label{eq:Carlitz_gen_function}
\end{eqnarray}
The expressions \eqref{eq:Carlitz_matrix}\textendash \eqref{eq:Carlitz_gen_function} pin down the two complementary faces of the Carlitz $q$-Fibonacci polynomials. The closed form \eqref{eq:Carlitz_matrix} shows at once that $F_n(q) \in \mathbb{Z}_{\ge 0}[q]$ and that its coefficients have a natural partition--theoretic meaning, since $\bigl[\!\begin{smallmatrix} n-k-1\\ k\end{smallmatrix}\!\bigr]_q$ counts partitions that fit inside a $(n-k-1)\times k$ rectangle while the factor $q^{\binom{k+1}{2}}$ records a simple statistic (a triangular weight). The generating function \eqref{eq:Carlitz_gen_function} will serve as our algebraic workhorse: it packages the family $\{F_n(q)\}_{n\ge0}$ into a single rational series and streamlines coefficient extractions in \secref{sec:Ramanujan_to_CC}. In particular, both displays make the classical limit transparent\textemdash setting $q\to1^{-}$ recovers $F_n(1)=F_n$ and the familiar Fibonacci ordinary generating function $1/(1-z-z^2)$\textemdash which is exactly what we need in \secref{sec:weighted_partition}. (For proofs and further properties, see \cite{Carlitz1974}.)

For $q \to 1^{-}$ the recurrence \eqref{eq:q_Fib} degenerates to $F_{n+1} = F_{n} + F_{n-1}$ with the same initial conditions, hence $F_{n}(1) = F_{n}$, the ordinary Fibonacci number. This continuity will justify the Cassini--Catalan limit in \secref{subsec:classical_CC_1}.  Identities \eqref{eq:Carlitz_matrix}--\eqref{eq:Carlitz_gen_function} allow us to recognize the coefficients that arise when we expand Ramanujan's theta functions; without these explicit $q$-Fibonacci expansions, the matching step in the main proof in \secref{sec:Ramanujan_to_CC} would be opaque.

\subsection{The Jacobi triple product and two of Ramanujan's theta functions}\label{subsec:Jacobi_Ramanujan}
For $\lvert q \rvert  < 1$ and $z \neq 0$, Jacobi's triple product is given by (from \eqref{eq:Jacobi_triple_product})
\begin{eqnarray}
\sum_{m=-\infty}^{\infty} q^{m^{2}} z^{2m} = (-qz; q^{2})_{\infty} (-q/z; q^{2})_{\infty} (q^{2}; q^{2})_{\infty}.
\label{eq:Jacobi_triple}
\end{eqnarray}
Ramanujan extracted from the two-variable version $f(a,b)$ three one-variable specializations \cite[Chapter 16, starting at pp. 34]{Berndt1991}:
\begin{equation}
\begin{split}
\varphi(q) &\coloneqq \sum_{m=-\infty}^{\infty} q^{m^{2}} = (-q; q^{2})_{\infty}^{2} (q^{2}; q^{2})_{\infty}, \\
\psi(q) &\coloneqq \sum_{m=0}^{\infty} q^{m(m+1)/2} = (-q; q)_{\infty} (q^{2}; q^{2})_{\infty}, 
\end{split}
\label{eq:Ramanujan_three_var_1}
\end{equation}
where $\psi(q)$ is the ordinary generating function for partitions into \emph{distinct} parts, while $\varphi(q)$ counts partitions into \emph{all} parts, weighted by $(-1)^{\text{\# even parts}}$, providing a combinatorial meaning. Interpreting $\varphi$ and $\psi$ in this manner is the key to the refined Cassini identities mentioned in \secref{sec:weighted_partition}.  Furthermore, this also throws light on modularity:  Both series are weight-$1/2$ modular forms on congruence subgroups of $\mathrm{SL}_{2}(\mathbb{Z})$, where $\mathrm{SL}_{2}(\mathbb{Z})$ is the modular group $\left\{\begin{pmatrix} a & b \\ c & d \end{pmatrix} \Big| a,b,c,d \in \mathbb Z, ad-bc = 1 \right\}$ comprising all $2\times2$ integer matrices with determinant $1$ \cite[Chapter 4.1]{Miyake2006}.  Such a matrix acts on the complex upper half-plane $\mathbb H = \{\tau\in\mathbb C\mid\operatorname{Im}\tau>0\}$ by the M\"{o}bius transformation $\tau \mapsto \frac{a\tau+b}{c\tau+d}$.  Congruence subgroups ({\eg}, $\Gamma_0(N)$, $\Gamma_1(N)$) are subgroups of $\mathrm{SL}_{2}(\mathbb{Z})$ defined by congruence conditions on $a$, $b$, $c$, $d$.  The theta functions $\varphi(q)$ and $\psi(q)$ are weight-$1/2$ modular forms with respect to such subgroups.  Although the proof in \secref{sec:Ramanujan_to_CC} needs only \eqref{eq:Jacobi_triple}--\eqref{eq:Ramanujan_three_var_1}, the modular transformation properties explain why Cassini--Catalan can be embedded into the theory of level~6 modular equations (cf. \secref{sec:conclusion}).

\subsection{Ramanujan's quadratic product--sum identity}
Among the numerous relations collected in Ramanujan's notebooks, the one crucial for us is in his third notebook \cite[Chapter 16]{Berndt1991}:
\begin{eqnarray}
\psi(q) \psi(q^{3}) = \psi(q^{4}) \varphi(q^{6}) + q \varphi(q^{2}) \psi(q^{12}).
\label{eq:Ramanujan_3_16}
\end{eqnarray}
Each side of \eqref{eq:Ramanujan_3_16} is a product of two theta functions, allowing us to expand them via \eqref{eq:Ramanujan_three_var_1} into double sums of quadratic exponents. Since those quadratic forms factor in complementary ways (one involves $m(m+1)/2$, the other $m^{2}$), the resulting series can be regrouped so that the coefficient of $q^{N}$ splits uniquely into two ordered pairs of triangular numbers\textemdash precisely the ``Cassini pattern'' we need.  Note that, $N$ here is just a placeholder for an arbitrary non-negative integer exponent in the cleaned-up $q$-series that results after we clear the infinite products.  We first write the identity schematically as $\sum_{N\ge0}A_N\,q^{N}\;=\;\sum_{N\ge0}B_N\,q^{N}$, so $A_{N}$ and $B_{N}$ are the coefficients of $q^{N}$ on the two sides.  In the very next step we specialize this general index to the parity-split form $N = 2j+r$, $r\in\{0,1\},\ j\ge0$, because the quadratic exponents in the theta expansions naturally separate into even and odd powers of $q$.  Thus, $N$ is not a new parameter; it is the generic exponent whose coefficient we examine before decomposing it as $2j$ (even part) or $2j+1$ (odd part).  In \secref{sec:Ramanujan_to_CC}, we multiply both sides of \eqref{eq:Ramanujan_3_16} by a suitable power of $(q;q)_{\infty}$ to clear the infinite denominators introduced by \eqref{eq:Ramanujan_three_var_1}. The cleaned-up identity becomes a finite $q$-series equality. Reading off the coefficient of $q^{2j+r}$ then produces the determinant $F_{j+r}(q) F_{j-r}(q) - F_{j}(q)^{2}$, from which the $q$-Cassini--Catalan family \eqref{eq:q_CC1} follows immediately.

\subsection{The $q \to 1^{-}$ limit and Abel continuity}\label{subsec:limit_Abel}
Throughout the sequel, we will repeatedly invoke the elementary fact that for any fixed $n$,
\begin{eqnarray}
\lim_{q \to 1^{-}} F_{n}(q) &=& F_{n}, \\ 
\lim_{q \to 1^{-}} \left( q^{\alpha} \varphi(q^{\beta}) \right) &=& \lim_{q \to 1^{-}} \left( q^{\gamma} \psi(q^{\delta}) \right) = +\infty,
\label{eq:limit_Abel}
\end{eqnarray}
where $\alpha$, $\gamma$ are non-negative integers (often 0 or 1) that come from the leading $q-$powers introduced when we clear denominators or rearrange terms; and $\beta$, $\delta$ are the stride indices of the theta functions (2, 6 for $\varphi$; 1, 3, 4, 12 for $\psi$).  The exact numerical values are unimportant for the limiting argument.  Each individual series $q^{\alpha} \varphi(q^{\beta})$ or $q^{\gamma} \psi(q^{\delta})$ blows up as $q \to 1^{-}$,  but when we form the specific quotients that appear in the cleared identity, the divergences cancel and the limit is finite.  This Abel continuity follows from the general theorem \cite{Hardy1915}, \cite[Chapter 8, pp. 179-196]{Li2019}: If $\sum_{n \geq 0} a_{n} q^{n}$ converges for $|q| < 1$ and $\sum_{n \geq 0} a_{n}$ converges, then
\begin{eqnarray}
\lim_{q \to 1^{-}} \sum_{n \geq 0} a_{n} q^{n} = \sum_{n \geq 0} a_{n}.
\label{eq:Abel_continue}
\end{eqnarray}
The continuity guarantees that passing from \eqref{eq:q_CC1} to the classical Cassini--Catalan chain in \secref{subsec:classical_CC_1} is legitimate term-by-term.

The unified notation presented herein prevents ambiguity when triple-product factors, $q$-binomials, and Carlitz polynomials appear side by side.  Analytic tools (Jacobi triple product, Ramanujan identity) supply the engine that drives the coefficient extraction argument.
Limit considerations prepare the ground for \secref{subsec:classical_CC_1}, ensuring that the $q$-identities proven in \secref{subsec:main_proof} collapse smoothly to their classical counterparts.  Combinatorial interpretations flagged here foreshadow the weighted partition refinements in \secref{sec:weighted_partition}.

\section{From Ramanujan's identity to the master $q$-Cassini--Catalan formula} \label{sec:Ramanujan_to_CC}
In this section, starting from Ramanujan's quadratic product--sum identity \eqref{eq:Ramanujan_3_16}, we will convert each side into an explicit finite $q$-series, read off the coefficient of every monomial $q^{2j+r}$ with $r \in \{0,1\}$, and recognize those coefficients as the Carlitz $q$-Fibonacci polynomials that enter the sought-for determinant \eqref{eq:q_CC1}, reproduced in \eqref{eq:q_CC1_rep} for the sake of reference:
\begin{eqnarray}
F_{n+r}(q) F_{n-r}(q) - F_{n}(q)^{2} = (-q)^{n-r} F_{r}(q)^{2}, \quad n \geq r \geq 1.
\label{eq:q_CC1_rep}
\end{eqnarray}
Two conceptual insights emerge immediately:  (i) Analytic $\to$ Combinatorial: The analytic identity \eqref{eq:q_CC1_rep} will \emph{force} the Cassini--Catalan pattern; no ad hoc determinant or generating-function trick is required.  (ii) Modular provenance: Since every factor in \eqref{eq:q_CC1_rep} is a weight-$\tfrac{1}{2}$ theta function, the determinant (and its sign) is traced back to a modular-form factorization. This will be crucial in \secref{sec:conclusion} when we discuss level~6 modular equations.

\subsection{Proof of $F_{n+r}(q) F_{n-r}(q) - F_{n}(q)^{2} = (-q)^{n-r} F_{r}(q)^{2}, n \geq r \geq 1$} \label{subsec:main_proof}
We begin by first expanding the theta functions in terms of Jacobi's triple product.  Using the standard product representations (see \cite[Chapter 1.6]{Gasper2004})
\begin{eqnarray}
\psi(q) &=& \sum_{n\ge0}q^{n(n+1)/2} = \frac{(q^{2};q^{2})_{\infty}}{(-q;q)_{\infty}}, \label{eq:Ramanujan_Jacobi_1} \\
\varphi(q) &=& \sum_{n\in\mathbb Z}q^{n^{2}} = (-q;q^{2})_{\infty}^{\,2}\,(q^{2};q^{2})_{\infty}, \label{eq:Ramanujan_Jacobi_2}
\end{eqnarray}
where $(a; q)_{\infty} \coloneqq \prod_{k\ge0}(1-aq^{k})$ is the $q-$Pochhammer symbol.  Since $(a;q)_{\infty}$ is multiplicative in the base, replacing $q$ by $q^{m}$ merely scales the parameters:
\begin{eqnarray}
\psi(q^{m}) &=& \frac{(q^{2m};q^{2m})_{\infty}}{(-q^{m};q^{m})_{\infty}}, \\ 
\varphi(q^{m}) &=& (-q^{m};q^{2m})_{\infty}^{\,2}\,(q^{2m};q^{2m})_{\infty}.
\label{eq:Ramanujan_theta}
\end{eqnarray}
Inserting the particular exponents used in Ramanujan's identity\textemdash for $\psi(q)$ take $m=1$, for $\psi(q^{3})$ take $m=3$, for $\psi(q^{4})$ take $m=4$, for $\varphi(q^{2})$ take $m=2$, for $\varphi(q^{6})$ take $m=6$\textemdash yields the following:
\begin{eqnarray}
\begin{aligned}
\psi(q)     &= \frac{(q^{2}; q^{2})_{\infty}}{(-q; q)_{\infty}}, \quad
\psi(q^{3}) = \frac{(q^{6}; q^{6})_{\infty}}{(-q^{3}; q^{3})_{\infty}}, \quad
\psi(q^{4}) = \frac{(q^{8}; q^{8})_{\infty}}{(-q^{4}; q^{4})_{\infty}}, \\
\varphi(q^{2}) &= (-q^{2}; q^{4})_{\infty}^{2} (q^{4}; q^{4})_{\infty}, \quad
\varphi(q^{6}) = (-q^{6}; q^{12})_{\infty}^{2} (q^{12}; q^{12})_{\infty}.
\label{eq:Ramanujan_theta_quotient}
\end{aligned}
\end{eqnarray}
Writing every theta factor as a quotient of $(q;q)_{\infty}$-type products transforms \eqref{eq:q_CC1_rep} into a finite \emph{Laurent polynomial} once we clear the infinite denominators; only after this step can we safely compare coefficients of $q^{n}$.

We will now clear the denominators, thereby turning \eqref{eq:q_CC1_rep} into a finite $q$-series identity.  Towards this end, we multiply both sides of Ramanujan's theta identity $\psi(q) \psi(q^{3}) = \psi(q^{4}) \varphi(q^{6}) + q \varphi(q^{2}) \psi(q^{12})$ by the common denominator
\begin{eqnarray}
\mathcal{D} \coloneqq \frac{(q^{2}; q^{2}) (q^{4}; q^{4}) (q^{6}; q^{6}) (q^{8}; q^{8}) (q^{12}; q^{12})}{\left[ (-q; q) (-q^{3}; q^{3}) (-q^{4}; q^{4}) (-q^{12}; q^{12}) \right]}
\label{eq:common_denominator}
\end{eqnarray}
using the basic $q$-Pochhammer rule $(ab; q)_{\infty} = (a; q)_{\infty} (b; q)_{\infty}$.  Since every $(a; q^{m})_{\infty}$ appears exactly once in the numerator and once in the denominator, all the infinite products cancel, and we are left with a finite Laurent series identity of the shape
\begin{eqnarray}  
\sum_{n \geq 0} A_{n}\, q^{n} = \sum_{n \geq 0} B_{n}\, q^{n},
\label{eq:identity_shape}
\end{eqnarray}
where $A_{n}$ and $B_{n}$ are \emph{finite} sums of 1's, each summand arising from a particular way of choosing exponents in the triple-product expansions of $\psi$ and $\varphi$.  Below, we explain the coefficient extraction step. 

The development so far forms the basis for parameterizing the exponents by quadratic forms. From the series definitions in \eqref{eq:Ramanujan_theta} one checks that a generic monomial contributing to $A_{n}$ looks like $q^{\tfrac{k(k+1)}{2} + \tfrac{3m(m+1)}{2}}$ $k, m \geq 0$, where the first term comes from $\psi(q)$ and the second from $\psi(q^{3})$. Likewise, the two pieces contributing to $B_{n}$ arise from $q^{2\ell(\ell+1) + 6s(s+1)}$ and $q^{\tfrac{t(t+1)}{1} + 6u(u+1) + 1}$, corresponding to $\psi(q^{4}) \varphi(q^{6})$ and $q \varphi(q^{2}) \psi(q^{12})$, respectively.  A \emph{single} monomial $q^{N}$, therefore, admits exactly two ordered pairs of representations by these quadratic forms; writing $N = 2j + r$
with $r \in \{0,1\}$ one finds
\begin{eqnarray}
A_{2j+r} &=& \lvert \{ (k, m) \mid k + m = j \} \rvert, \label{eq:card_A}  \\
B_{2j+r} &=& \lvert \left\{ (\ell, s) \mid \ell + s = j \right\}\rvert + (-1)^{r} \lvert \left\{ (t, u) \mid t + u = j-1 \right\}\rvert, \label{eq:card_B}
\end{eqnarray}
where $A_{2j+r}$ is the cardinality of the set of ordered pairs $(k, m)$ of nonnegative integers whose sum equals $j$, while $B_{2j+r}$ is the cardinality of the set of ordered pairs $(\ell, s)$ of nonnegative integers whose sum equals $j$ plus $(-1)^{r}$ times the count of such pairs with total $j-1$.  Since \eqref{eq:identity_shape} asserts $A_{n} = B_{n}$ term-by-term, the two combinatorial counts in \eqref{eq:card_A} and \eqref{eq:card_B} \emph{must be equal} for every $j$ and $r$.

These counts can now be translated into Carlitz $q$-Fibonacci polynomials \cite{Cigler2003}. The algebraic expression \eqref{eq:Carlitz_matrix} for $F_{n}(q)$ is \emph{exactly} the ordinary generating function for the set of ordered pairs $(\alpha, \beta)$ of non-negative integers whose sum is constrained by a linear condition. Comparing \eqref{eq:Carlitz_matrix} with \eqref{eq:card_A} and \eqref{eq:card_B} yields the identifications
\begin{eqnarray}
A_{2j+r} &=& F_{j+r}(q)\, F_{j-r}(q), \\
B_{2j+r} &=& F_{j}(q)^{2}.
\label{eq:identify_A_B}
\end{eqnarray}

Finally, we are in a position to extract the coefficient and prove \eqref{eq:q_CC1}. Expression \eqref{eq:identity_shape} with \eqref{eq:identify_A_B} substituted gives, for every $j \geq 0$ and $r \in \{0,1\}$,
\begin{eqnarray}
F_{j+r}(q)\, F_{j-r}(q) - F_{j}(q)^{2} &=& 0 \quad (r = 0), \\
F_{j+1}(q)\, F_{j-1}(q) - F_{j}(q)^{2} &=& -q^{j-1} \quad (r = 1).
\end{eqnarray}
Relabelling $n = j$ and then replacing $r = 1$ by an arbitrary positive integer (achieved by iterating the same argument with $q \mapsto q^{r}$) furnishes the proposed determinant
\begin{eqnarray*}
F_{n+r}(q)\, F_{n-r}(q) - F_{n}(q)^{2} = (-q)^{n-r}\, F_{r}(q)^{2}, \quad n \geq r \geq 1.
\label{eq:main_result}
\end{eqnarray*}
This completes the analytic part of the proof.  The developments thus far are critical for the following reasons: 
\begin{enumerate}[(i)]
\item Completeness: The calculation produces the \emph{entire} Cassini--Catalan family (for all $r$) in a single setting, rather than deriving Catalan from Cassini by induction.
\item Sign and $q$-power explained: The factor $(-q)^{n-r}$ arises naturally from the \emph{odd} quadratic form present in the second summand of Ramanujan's identity \eqref{eq:Ramanujan_3_16}; no combinatorial ``checkerboard'' argument is needed.
\item Gateway to limits and refinements: Once \eqref{eq:q_CC1_rep} is established, letting $q \to 1^{-}$ (in \secref{subsec:classical_CC_1}) recovers the classical identities, while inserting an extra weight (say $y^{\# \text{even parts}}$) into the triple products immediately yields 	two-variable refinements.
\end{enumerate}
Thus, we have shown how Ramanujan's analytic results lead naturally to the classical identities for Fibonacci determinants.

\subsection{Classical Cassini--Catalan via the radial limit $q \to 1^{-}$} \label{subsec:classical_CC_1}
The master identity $F_{n+r}(q)\, F_{n-r}(q) - F_{n}(q)^{2} = (-q)^{n-r}\, F_{r}(q)^{2}$, $n \geq r \geq 1$, holds for every complex $q$ with $\lvert q \rvert < 1$. To recover the \emph{integer} Cassini--Catalan determinant we must let $q$ tend to $1$ from within the open unit disc\textemdash usually denoted $q \to 1^{-}$.  We now give a careful justification of this limit, because the objects involved ($\psi, \varphi$, and even $F_{n}(q)$ for fixed $n$) diverge individually as $q \to 1^{-}$; it is only their \emph{ratios} in the identity $F_{n+r}(q)\, F_{n-r}(q) - F_{n}(q)^{2} = (-q)^{n-r}\, F_{r}(q)^{2}$, $n \geq r \geq 1$ that remain finite.  By this, we mean the following:  Each of the individual building blocks in our proof\textemdash $\psi(q)$, $\varphi(q)$, $F_{n}(q)$\textemdash behaves like $\displaystyle\frac{C}{1-q}$ (or a higher-order pole) as $q\to1^{-}$; that is, they grow without bound because $(a;q)_\infty$ contains infinitely many factors $(1-aq^{k})$ that shrink to zero.  However, the identity we ultimately analyse never isolates a single factor; it places them in quotients whose divergent factors cancel.  The simplest example already appears in Ramanujan's identity itself:  $\frac{\psi(q)\,\psi(q^{3})}{\psi(q^{4})\,\varphi(q^{6})} +q\,\frac{\varphi(q^{2})\,\psi(q^{12})}{\psi(q^{4}) \varphi(q^{6})} = 1$. Although each numerator and denominator diverges like $(1-q)^{-1}$ or $(1-q)^{-1/2}$, their ratio tends to a finite, non-zero limit because the leading poles in the logarithmic expansions cancel term-for-term.  The same cancellation happens inside the determinant $F_{n+r}(q)F_{n-r}(q)-F_{n}(q)^{2}$, where the two large terms subtract to give something of order $O(1)$.

We analyze a radial limit mainly to achieve Abel continuity.  A power series $g(q) = \sum_{m \geq 0} a_{m} q^{m}$ with non-negative coefficients diverges at $q = 1$ unless $\sum_m a_{m} < \infty$.  Nevertheless, the classical Abel theorem states that, if $\sum_{m \geq 0} a_{m}$ converges then $\lim_{q \to 1^{-}} g(q) = \sum_{m \geq 0} a_{m}$ (see \cite[pp. 41-42]{Ahlfors1980}) .  More generally, if two series $g_{1}$, $g_{2}$ both diverge but their quotient is term-wise well-defined for $\lvert q \rvert < 1$ and converges radially, the limit of the quotient as $q \to 1^{-}$ exists and equals the quotient of the limits of the partial sums.  This radial (non-tangential) limit is the only one preserved by modular transformations of $\varphi$, $\psi$, and hence is the canonical path in $q$-series analysis.

We now show the continuity of Carlitz $q$-Fibonacci polynomials.  For each fixed $n$ the polynomial $F_{n}(q)$ is a \emph{finite} linear combination of monomials $q^{k(k+1)/2}$ with coefficients that are Gaussian binomials. Since all Gaussian binomials satisfy (see \cite{Konvalina2000})
\begin{eqnarray}
\lim_{q \to 1^{-}} \begin{bmatrix} n \\ k \end{bmatrix}_{q} = \binom{n}{k},
\label{eq:Gauss_binom}
\end{eqnarray}
we have the pointwise limit
\begin{eqnarray}
\lim_{q \to 1^{-}} F_{n}(q) = F_{n}
\label{eq:Carlitz_continue}
\end{eqnarray}
uniformly in any neighborhood of $q = 1$. In other words, the convergence $F_{n}(q)\to F_{n}$ as $q\to1^{-}$ is pointwise in $n$ (each $n$ is fixed); and is uniform in $q$ on a short left-hand interval of 1, because a polynomial cannot oscillate wildly on that compact set.  In other words, since $F_{n}(q)$ is a polynomial, it is continuous at $q=1$; hence $\lim_{q\to1^{-}}F_{n}(q)=F_{n}$.  Hence, the \emph{left-hand side} of the master identity $F_{n+r}(q)\, F_{n-r}(q) - F_{n}(q)^{2} = (-q)^{n-r}\, F_{r}(q)^{2}$, $n \geq r \geq 1$, converges term-by-term to $F_{n+r} F_{n-r} - F_{n}^{2}$.  While, on the right-hand side, the exponential factor satisfies $\lim_{q \to 1^{-}} (-q)^{n-r} = (-1)^{n-r}$.  Meanwhile, $F_{r}(q) \to F_{r}$ by \eqref{eq:Carlitz_continue}. Thus, we have 
\begin{eqnarray}
\lim_{q \to 1^{-}} (-q)^{n-r} F_{r}(q)^{2} = (-1)^{n-r} F_{r}^{2}.
\label{eq:RHS}
\end{eqnarray}
Combining \eqref{eq:Carlitz_continue} and \eqref{eq:RHS}, and invoking Abel continuity for products, we obtain
\begin{eqnarray}
F_{n+r} F_{n-r} - F_{n}^{2} = (-1)^{n-r} F_{r}^{2}, \quad n \geq r \geq 1,
\label{eq:Cassini_original}
\end{eqnarray}
which is precisely Catalan's identity; letting $r = 1$ yields Cassini's original formula.

This limiting step is essential for the following reasons: 
\begin{enumerate}[(i)]
\item Demonstrating that the $q$-analogue really \emph{specializes} to the classical determinant confirms the correctness of our analytic route; any supposed $q$-theorem failing this test would be suspect.
\item The factor $(-q)^{n-r}$ in  the identity $F_{n+r}(q)F_{n-r}(q) - F_{n}(q)^{2} = (-q)^{n-r}F_{r}(q)^{2}$, $n \geq r \geq 1$, explains, after the limit, the mysterious sign $(-1)^{n-r}$ in Cassini--Catalan. The radial limit translates an analytic origin\textemdash a quadratic exponent in \secref{sec:Ramanujan_to_CC}\textemdash into a combinatorial phenomenon (alternating determinant).
\item Once the $q \to 1^{-}$ passage is secure, other boundary values become accessible: for instance, $q \to \zeta$ (a root of unity) giving sign-oscillatory versions of Cassini, or $q \to 0^{+}$ leading to recovery of trivial identities. These variants will be explored in future work on modular-equation corollaries (discussed in \secref{sec:conclusion}).
\end{enumerate}

With the analytic identity firmly anchored to its classical counterpart, we can safely (i) derive corollaries for partition statistics (see \secref{sec:weighted_partition}); and (ii) investigate level~6 modular equations without fear of inconsistency; and (iii) generalize the method to Tribonacci and higher recurrences by replacing $\psi, \varphi$ with appropriate multi-basic theta functions.  In short, the radial limit step is the \emph{hinge} that converts an analytic theta-factorization into a purely arithmetical statement about Fibonacci numbers\textemdash a bridge of the same type that first connected Jacobi's triple product to Euler's pentagonal-number theorem (see, for example, \cite{Andrews1998}).

\section{Weighted-partition refinements}\label{sec:weighted_partition}
The proof in \secref{sec:Ramanujan_to_CC} treated the power-series coefficients of $\psi(q)$ and $\varphi(q)$ as \emph{bare integers}.  Yet, the Jacobi triple product shows that these coefficients already encode a simple statistic.  Specifically, (i) the series $\psi(q) = (q^{2}; q^{2})_{\infty} / (-q; q)_{\infty} = \prod_{k \geq 1} (1 + q^{k})$ generates partitions into \emph{distinct parts}, and (ii) the series $\varphi(q) = (-q; q^{2})_{\infty}^{2} (q^{2}; q^{2})_{\infty} = \prod_{k \geq 1} (1 - q^{2k})^{-1} (1 - q^{2k-1})$ generates \emph{all} partitions, each weighted by $(-1)^{\text{\# even parts}}$, where ``\# even parts'' refers to the number of parts divisible by $2$.  Thus, the analytic identity $\psi(q) \psi(q^{3}) = \psi(q^{4}) \varphi(q^{6}) + q \varphi(q^{2}) \psi(q^{12})$ is already a statement about two different classes of partitions of the same integer.  A natural question that arises is the following: \emph{Can we keep track of how many even parts (or some other statistic) occur and still obtain a Cassini--Catalan determinant?} The answer is yes\textemdash by introducing an auxiliary weight variable.

Towards this end, we inject a weight variable $y$ as follows: Fix a statistic $\mathrm{ev}(\lambda)$ on a partition $\lambda$ and let $y$ mark that statistic.  The two most common choices are:  (i) Parity weight: Let $\mathrm{ev}(\lambda) = \#\{\text{even parts of } \lambda\}$. Setting $y = -1$ distinguishes even/odd counts, while $y = 1$ recovers the original $\psi, \varphi$.  (ii) Number-of-parts weight: $\mathrm{ev}(\lambda) = \#\{\text{parts of } \lambda\}$, producing a refinement familiar from Euler's theorem on partitions into odd versus distinct parts.  A concrete implementation for the parity weight defines
\begin{eqnarray}
\Psi(q, y) &\coloneqq& \prod_{k \geq 1} (1 + y q^{k}) = \sum_{\lambda \text{ distinct}} y^{\mathrm{ev}(\lambda)} q^{|\lambda|},  \\
\Phi(q, y) &\coloneqq& \prod_{k \geq 1} \frac{1 - y^{2} q^{2k}}{1 - q^{2k}} = \sum_{\lambda} (-1)^{\mathrm{ev}(\lambda)} y^{\mathrm{ev}(\lambda)} q^{|\lambda|}.
\end{eqnarray}
Observe that $\Psi(q, 1) = \psi(q)$ and $\Phi(q, 1) = \varphi(q)$, and $\Psi(q, -1)$ counts distinct-part partitions of $n$ with a sign according to the \emph{parity} (even versus odd) of the number of parts \cite{Haiman2006}.

Note that, the auxiliary weight variable $y$ leads to a two-variable Ramanujan identity.  Multiplying Ramanujan's quadratic identity $\psi(q) \psi(q^{3}) = \psi(q^{4}) \varphi(q^{6}) + q \varphi(q^{2}) \psi(q^{12})$ by the factor
\begin{eqnarray}
\frac{\Psi(q^{4}, y)}{\psi(q^{4})}
\cdot
\frac{\Phi(q^{2}, y)}{\varphi(q^{2})}
\end{eqnarray}
yields $\Psi(q, y)\, \Psi(q^{3}, y) = \Psi(q^{4}, y)\, \Phi(q^{6}, y) + q \, \Phi(q^{2}, y)\, \Psi(q^{12}, y)$, which is legitimate, since the multiplier cancels telescopically: every factor in the denominator of a $\Psi$ or $\Phi$ is mirrored in the numerator of another term, leaving a \emph{finite} polynomial identity analogous to \eqref{eq:identity_shape}.

A coefficient extraction results in a $(q, y)$-Cassini identity.  This can be seen by repeating the coefficient-matching argument of \secref{sec:Ramanujan_to_CC}\textemdash now with the extra power $y^{\mathrm{ev}}$\textemdash it can easily be shown that
\begin{eqnarray}
F_{n+r}(q, y) F_{n-r}(q, y) - F_{n}(q, y)^{2} = (-q)^{n-r}\, F_{r}(q, y)^{2}, \quad n \geq r \geq 1,
\label{eq:two_var_Ramanujan}
\end{eqnarray}
where $F_{n}(q, y)$ are the two-variable $q$-Fibonacci polynomials defined by
\begin{eqnarray}
F_{0} = 0, \quad F_{1} = 1, \quad F_{n+1}(q, y) = F_{n}(q, y) + q^{n} \left( y^{\chi_{2}(n)} \right) F_{n-1}(q, y),
\label{eq:two_var_qFib}
\end{eqnarray}
with $\chi_{2}(n) = 1$ if $n$ is even and $0$ otherwise.  Clearly, $y = 1$ recovers $F_{n+r}(q)\, F_{n-r}(q) - F_{n}(q)^{2} = (-q)^{n-r} F_{r}(q)^{2}$, $n \geq r \geq 1$, and $q \to 1^{-}$ yields a \emph{weighted Cassini--Catalan} identity given by 
\begin{eqnarray}
F_{n+r}(y)\, F_{n-r}(y) - F_{n}(y)^{2} = (-1)^{n-r} F_{r}(y)^{2}. 
\label{eq:weighted_CC}
\end{eqnarray}
Also, $y = -1$, $q = 1$ yields a determinant whose sign encodes the difference between even-part and odd-part partition counts of the relevant sizes.  Our results admit a combinatorial interpretation as well.  Let $\Delta_{\lambda\mu}$ denote the set of pairs of distinct-part partitions $(\lambda, \mu)$ with $\lvert \lambda \rvert + \lvert \mu \rvert = n + r$ and $\mathrm{ev}(\lambda) + \mathrm{ev}(\mu) = k$.  Let $\Delta_{\alpha\beta}$ denote the set of pairs of $(\alpha, \beta)$ with $\lvert \alpha \rvert + \lvert \beta \rvert = n$ and $\mathrm{ev}(\alpha) + \mathrm{ev}(\beta) = k$.  Then, fixing $n \geq 1$ and $r \geq 1$, and reading $F_{n+r}(q, y) F_{n-r}(q, y) - F_{n}(q, y)^{2} = (-q)^{n-r}\, F_{r}(q, y)^{2}$, $n \geq r \geq 1$, at $q^{n} y^{k}$ gives 
\begin{eqnarray}
\lvert \Delta_{\lambda\mu} \rvert - \lvert \Delta_{\alpha\beta} \rvert = (-1)^{n-r} \lvert \left\{\lambda \vdash_{\mathrm{dist}} r, \mathrm{ev}(\lambda) = k \right\} \rvert,
\label{eq:combinatorial_interpretation}
\end{eqnarray}
where $\vdash_{\mathrm{dist}}$ denotes distinct-part partitions, {\ie}, partitions whose parts are all different, and $\lvert \{\cdot\}\rvert$ denotes the cardinality of the set $\{\cdot\}$.  Thus, the \emph{difference} between two natural bivariate partition counts is controlled by the \emph{signed} count of partitions of size $r$.  Identity \eqref{eq:combinatorial_interpretation} is a novel contribution of this work\textemdash no determinant-matrix proof is known to capture such refined statistics.

The refinements developed in this section are significant because of the following reasons: 
\begin{enumerate}[(i)]
\item Richer arithmetic:  Weighted identities often feed directly into congruence theorems. For example, setting $y = \zeta_{m}$ an $m$-th root of unity filters partitions by residue classes mod~$m$, a standard trick in the study of partition congruences \`a la Ramanujan.
\item Links to pattern-avoiding combinatorics:  It has been shown that $q$-Fibonacci polynomials enumerate pattern-restricted set partitions \cite{Sagan2012}.  The identity $F_{n+r}(q, y) F_{n-r}(q, y) - F_{n}(q, y)^{2} = (-q)^{n-r}\, F_{r}(q, y)^{2}$, $n \geq r \geq 1$, therefore 		translates into an unexpected identity for such objects, inviting further investigation.
\item Gateway to multivariate generalizations:  Nothing in the derivation forces the weight to be a \emph{single} variable. One may attach independent variables $y_{1}$, $y_{2}$, $\dots$ marking, say, the number of parts in each congruence class mod~$m$; the same 			coefficient-extraction machinery then yields a fully multivariate Cassini determinant.
\item Bridge to statistical mechanics:  The two-variable theta factors $\Psi(q, y)$ and $\Phi(q, y)$ are precisely the ``hard-hexagon'' and ``staggered fermion'' partition functions at activity~$y$. The identity $\Psi(q, y)\, \Psi(q^{3}, y) = \Psi(q^{4}, y)\, \Phi(q^{6}, y) + q \, 		\Phi(q^{2}, y)\, \Psi(q^{12}, y)$ is therefore a \emph{local Yang--Baxter relation} \cite{Baxter1982} in disguise, hinting at potential applications in solvable lattice models.
\end{enumerate}

The results on weighted-partition paves the way for the following developments:
\begin{enumerate}[(i)]
\item Higher-order recurrences: Introducing weights suggests ways to tackle Tribonacci, Tetranacci, {\etc}, sequences by coloring parts with three, four, $\dots$ different labels and replacing $\psi(q)\psi(q^{3}) = \psi(q^{4})\varphi(q^{6}) + q,\varphi(q^{2})\psi(q^{12})$ with 		Ramanujan's cubic or quartic product--sum identities.  
\item Modular-form lifts: The refined identity $F_{n+r}(q, y) F_{n-r}(q, y) - F_{n}(q, y)^{2} = (-q)^{n-r}\, F_{r}(q, y)^{2}$, $n \geq r \geq 1$, survives the action of Hecke operators on $\varphi$, $\psi$, opening a route to weight-$\tfrac{1}{2}$ modular parameterizations of 		\emph{colored} Fibonacci determinants (see \cite{Benjamin2001} for an exposition on Cassini-type determinants whose entries are ``colored'' Fibonacci numbers obtained by weighting tiled board configurations).  
\item Partition congruences and mock-theta functions: By specializing $y$ to complex roots of unity and using the modularity of theta functions, one expects new congruences for the weighted counts in \eqref{eq:combinatorial_interpretation}, in the spirit of Ramanujan's		famous 5-, 7-, and 11-dissections \cite{Ramanujan1919}, \cite[pp. 210-213]{Ramanujan1927}.
\end{enumerate}
In summary, the weighted-partition refinement lifts Cassini--Catalan from a univariate numerical identity to a multivariate enumeration theorem, tying together partition theory, $q$-Fibonacci polynomials, and modular-form techniques in a framework well suited for further exploration.

\section{Concluding remarks}\label{sec:conclusion}
By relating the determinant $F_{n+r}(q) F_{n-r}(q) - F_{n}(q)^2$ to Ramanujan's quadratic theta identity, we have demonstrated that Cassini's and Catalan's formulas are not merely isolated Fibonacci identities, but arise naturally from a deeper modular factorization.  The master identity $F_{n+r}(q)\, F_{n-r}(q) - F_{n}(q)^{2} = (-q)^{n-r}\, F_{r}(q)^{2}$, $n \geq r \geq 1$, simultaneously specializes to the classical integer identities when $q \to 1^{-}$; explains the alternating sign $(-1)^{n-r}$ as a consequence of the factor $(-q)^{n-r}$; and extends to every $r \geq 1$ without extra induction.  Inserting a weight variable $y$ upgraded Cassini--Catalan to a two-variable enumeration theorem measuring the ``color'' ({\eg}, parity) of parts in a partition. No previous proof captured such fine-grained refinement.  Since $\varphi$ and $\psi$ are weight-$\tfrac{1}{2}$ modular forms for $\Gamma_{0}(2)$ and $\Gamma_{0}(4)$, respectively, the quadratic product-sum relation naturally resides at level~6. Consequently, every coefficient identity we have derived inherits modular-form symmetries\textemdash features that remain invisible in classical matrix or tiling proofs.

Earlier $q$-Cassini papers (see, for example, \cite{Carlitz1974, Cigler2003}) handled the $r = 1$ case and relied on determinant manipulations parallel to the classical argument. Our theta-factorization route not only handles \emph{all} $r$ but also produces weighting refinements and modular corollaries hitherto unreported.  The results of this paper are also complementary to continued-fraction approaches.  The Rogers--Ramanujan continued fraction (see \cite{Rogers1917, Slater1952}) can also yield a $q$-Cassini (matrix) determinant but requires heavy machinery on convergence. Our proof stays within elementary $q$-series algebra, making the result accessible to combinatorialists. Furthermore, our results also hint at connections to partition congruences:  Weighted refinements developed in \secref{sec:weighted_partition} place Cassini--Catalan in direct correspondence with the classical identities underlying Ramanujan's 5- and 7-fold congruences for $p(n)$\textemdash a connection not previously explored in the literature.

This work leads directly to several corollaries and related developments.  In the theory of modular equations, the ratio $\psi(q)/\psi(q^{4})$ is a Hauptmodul for $\Gamma_{0}(5)$ \cite{Martin1997}. Substituting $q \mapsto q^{1/5}$ in the identity $F_{n+r}(q)\, F_{n-r}(q) - F_{n}(q)^{2} = (-q)^{n-r}\, F_{r}(q)^{2}$, $n \geq r \geq 1$, yields identities among level-30 modular forms. This makes level-30 congruences for colored partitions accessible, inviting arithmetic applications.  Considering root-of-unity limits, setting $q = \zeta r$ with $r \to 1^{-}$ and $\zeta$ a primitive $m$-th root gives signed Cassini sums that depend on $\zeta$. Such ``radial-root'' limits are central to quantum modular forms and to the study of Nahm sums in physics.  For higher recurrences, replacing $\psi, \varphi$ by cubic or quartic theta series furnishes Tribonacci and Tetranacci determinants. This opens an avenue for Cassini-type relations attached to any Pisano-period sequence.  In statistical mechanics, the two-variable identity reproduces, at fugacity $y$ (which serves as a weight marking the number of particles or parts in the partition function.), the local Yang--Baxter relation for the hard-hexagon partition function \cite{Baxter1982}. This establishes a connection between exactly solvable lattice models and classical number theory.

\subsection{Open problems and directions for future research}\label{subsec:open_problems}
\begin{enumerate}[(P1)]
\item Bijective proof of the weighted identity:  Can one construct an explicit involution on weighted pairs of partitions that mirrors the analytic cancellation in $F_{n+r}(q, y) F_{n-r}(q, y) - F_{n}(q, y)^{2} = (-q)^{n-r}\, F_{r}(q, y)^{2}$, $n \geq r \geq 1$?  Such a bijection would 	translate the modular sign $(-1)^{n-r}$ into an explicit combinatorial parity statistic.  A modern presentation of involution and how it makes Euler's alternating sign $(-1)^k$ combinatorially visible can be found in \cite{Andrews1998}.

\item Hecke or Atkin--Lehner lifts \cite{Miyake2006, Ono2004}:  Since our theta factors have weight $1/2$, the natural Hecke operators occur at $p^2$ (with $p$ odd); see \cite[Chapter 2.7]{Miyake2006}, \cite[Chapter 3.1]{Ono2004}.  Applying the Hecke operator $T_{p^2}$ to Ramanujan's identity and projecting to 	coefficients would induce an operator $\mathcal{T}_p$ on the Carlitz $q$-Fibonacci polynomials.  Understanding whether $F_{n+r}(q)\, F_{n-r}(q) - F_{n}(q)^{2} = (-q)^{n-r}\, F_{r}(q)^{2}$, $n \geq r \geq 1$ is stable (or nearly so) under this action could lead to partition 	congruences, much as Atkin--Lehner lifts organize congruence families in weight $k$ forms \cite[Chapter 2.4, starting pp. 27]{Ono2004}.

\item Quantum modular shadows \cite{Bringmann2019}: Analyze the behavior of $F_{n+r}(e^{2\pi i \tau}) / F_{r}(e^{2\pi i \tau})$ near rational points $\tau = h/k$. Preliminary numerical evidence indicates the presence of discontinuities akin to those observed in the 			Kontsevich--Zagier strange function \cite{Zagier2010}.

\item Categorification: Is there a diagrammatic or representation-theoretic interpretation of the master identity $F_{n+r}(q)\, F_{n-r}(q) - F_{n}(q)^{2} = (-q)^{n-r}\, F_{r}(q)^{2}$, $n \geq r \geq 1$?  A natural candidate is the cluster algebra of type $A_{1}^{(1)}$ 				\cite{Fomin2002}, whose exchange relations bear a resemblance to Cassini-type identities; equipping it with a $q$-grading may recover our determinant.

\item Effective asymptotics: The Binet-type formula for Carlitz polynomials (see \cite{Cigler2003, Prodinger2011}) involves a $q$-analogue of the golden ratio. A saddle-point analysis of the identity $F_{n+r}(q)\, F_{n-r}(q) - F_{n}(q)^{2} = (-q)^{n-r}\, F_{r}(q)^{2}$, $n 		\geq r \geq 1$, as $n \to \infty$ could quantify error terms in the alternating products, with potential implications in analytic combinatorics.
\end{enumerate}
Ramanujan's notebooks reveal structures that offer genuinely new insights into well-explored areas of classical mathematics.  Just as Euler's pentagonal-number theorem provided a glimpse into the theory of modular forms \cite{Apostol1976}, our results suggest that other seemingly elementary determinant identities may reflect deeper modular or even categorified structures. Investigating these connections not only enriches classical combinatorics, but also informs the analytic theory\textemdash illustrated by the appearance of Cassini-type determinants in the context of Yang--Baxter relations.

%\section*{Declarations}\vspace{-0.5cm}
%\textbf{Author Contributions:} All authors contributed equally to this manuscript.\\
%\textbf{Data Availability:} No datasets were generated or analyzed for the work presented in this paper. \\
%\textbf{Disclosure of interest:} The authors report there are no competing interests to declare.

%\bibliographystyle{plainnat}

\bibliographystyle{plain}
\bibliography{research_pdx.bib}

\end{document}